\documentclass[ajp, preprint, author-numerical]{revtex4-1} 

\usepackage{amsmath} 
\usepackage{amsfonts} 

\usepackage{graphicx} 
\usepackage{threeparttable} 

\begin{document}
\def \RR {\mathbb{R}}
\def \GG {\mathbb{G}}

\title{Extensions of real numbers using coset groups}

\author{Horia I. Petrache}
 
\affiliation{Department of Physics, Indiana University Purdue University Indianapolis, Indianapolis, IN 46202}

\email{hpetrach@iupui.edu}

\today


\begin{abstract}

Extensions of real numbers in more than two dimensions, in particular quaternions and octonions are finding applications in physics due to the fact that they naturally capture certain symmetries of physical systems. Here it is shown that the property of closure of coset groups can be used to generate the basis and general multiplication rules for extensions of real numbers in a systematic way. The coset approach has the advantage that multiplication rules follow directly from group closure instead of being postulated. In this approach, constraints on multiplication parameters can be formulated in ways that capture the symmetry features of the coset group. A complete classification of numbers systems is therefore obtained based on possible group structures of a given order. General matrix representations are also obtained through the coset procedure  and by construction, the form of these matrices is invariant under matrix addition and multiplication. Since group symmetries are captured naturally into each number system, the coset group approach can add insight into the utility of multidimensional number systems in describing symmetries in nature.
\end{abstract}

\maketitle

\section{Introduction}

The simplest, most common example of an extension of real numbers are the complex numbers omnipresent in modern physics. Mathematically, complex numbers are constructed as an algebra over the field of reals using a set of two elements $1$ and $i$ (called a basis) and by postulating the multiplication rule $i^2 = -1$ \cite{yaglom, hungerford, kantor89}. The important feature is that the set of complex numbers is closed under both addition and multiplication. Other kinds of 2-dimensional extensions of real numbers can be constructed by modifying the multiplication rule. For example, split-complex numbers correspond to $i^2 = 1$ (also called hyperbolic numbers) and dual complex numbers correspond to $i^2 = 0$ \cite{yaglom, hucks}. 

	Another example gaining general interest in physics are the quaternions \cite{adler}. These are extensions of real numbers in 4 dimensions constructed using the basis $\{1, i, j, k\}$ with multiplication rules $i^2 = j^2 = k^2 = ijk = -1$ \cite{hamilton}. Just like the set of complex numbers, the set of quaternions is closed under both addition and multiplication. Due to their particular group structure, quaternions can be used in physics for representations of the Lorentz group and in general for transformations involving 4-vectors \cite{manogue-lorentz-transf, deleo}. Quaternions are also receiving attention in quantum mechanics due to a direct relationship with Pauli matrices and Pauli's group \cite{finkelstein, ISI:000265946900055, ISI:000262251100071, ISI:000293248600006}. Other variations of quaternions have been postulated, such as split-quaternions \cite{ISI:000275898500003} and bicomplex numbers which are obtained from a 4-dimensional basis like the quaternions but with different multiplication rules \cite{yaglom}. 

The emerging role of quaternions in physics suggests that other extensions of real numbers could possibly be useful for describing symmetries in nature \cite{conway-smith}. However, what extensions are possible and how to "discover" the appropriate multiplication rules? As shown here, a systematic approach can be based on the concept of coset products from group theory. Cosets are shifted replicas of a given subgroup forming a set that is closed under coset product (see Appendix for a brief introduction). It is this property of closure of the coset group that gives rise naturally to consistent multiplication rules for multidimensional numbers. In this respect, the coset approach is an alternative to the established algebra over fields theory with the advantage that multiplication rules follow directly from group closure instead of being postulated. A complete classification of multiplication rules is obtained in this way. Since all group structures are known (especially for small groups) multiplication rules are fully determined. 

In what follows, the general coset construction is presented and then applied to coset groups of order $n = 2$, $3$, and $4$. Multiplication rules and general matrix representations are derived for each number system corresponding to these small groups. By construction, the form of these matrices is invariant under matrix addition and multiplication. For additional reading on multidimensional number systems, references \cite{hawkins, kantor89, mcconnel98} review hypercomplex numbers in general, reference \cite{conway-smith} discusses quaternions and octonions, and references \cite{olariu-3d, olariu-4d, olariu-book} present a detailed analysis of a selection of multidimensional numbers including geometrical representations and functional analysis. For a historical perspective on the early development of quaternions, readers are directed to the original paper by Hamilton \cite{hamilton} and to references \cite{hardy-quater, hathaway, hawkes, taber, dickson-div-alg, dickson-quater,  mcconnel98}.

\section{General construction of coset extensions of reals}

Consider that there exists a set of elements $g_i$ outside of the set of real numbers $\RR$ but compatible with operations in $\RR$. Using these elements, we can construct a collection $\GG$ of cosest of $\RR$, 
\begin{equation}
\GG = \{...g_i\RR, g_ig_j\RR, g_ig_jg_k\RR, ...\},
\label{eq:cosets}
\end{equation}
where $g_i\RR$ represents the set obtained by multiplying each real number by $g_i$. In $g_i\RR$, all real numbers have been "shifted" or scaled by a factor of $g_i$, in $g_ig_j\RR$ by a factor of $g_ig_j$ and so forth. These shifted sets are called cosets (see Appendix). The set $\GG$ is then the set of cosets of $\RR$ and we require that $\GG$ is finite and closed under coset product, namely if $g_i\RR \in \GG$ and $g_j\RR \in \GG$, then $g_ig_j\RR \in \GG$ for all elements $g_i$ and $g_j$. (Technically, $\GG$ is a quotient group modulo $\RR$ of a larger group). 

By restricting the number $n$ of distinct elements in $\GG$, we generate various extensions of real numbers in which multiplication rules follow from the structure of the coset group. We will treat explicitly the simplest cases in which $\GG$ is one of the small groups with $n =$ 2, 3, or 4 elements. From the general theory of finite groups we know that for $n = 2$ and $n = 3$ the only possible group structure is cyclic, while for $n = 4$ there are two distinct groups, one cyclic and the other called the Klein four-group. By setting $\GG$ to be one of these groups, we obtain extensions of $\RR$ in a systematic way. Multiplication rules follow naturally from the closure of the coset group showing explicitly the underlying symmetry of the number system. 

\section{Case $n = 2$} 

In this simplest case we require that $\GG$ has only two distinct elements, $\GG = \{\RR, g\RR\}$ and is closed under coset product. As expected, this construction generates complex numbers. Although complex numbers can be defined and constructed in a number of ways as in abstract algebra textbooks (\textit{e. g.} references \cite{yaglom,hungerford}), the coset construction is presented here in detail in order to illustrate this particular procedure using a familiar number system.

Imposing that $\GG = \{\RR, g\RR\}$ is closed requires that each element in Eq.~\ref{eq:cosets} must be either $\RR$ or $g\RR$. Assuming that $g^{-1}\RR = \RR$ and multiplying by $g$, we obtain $gg^{-1}\RR = g\RR$ which gives $\RR = g\RR$. This contradicts the requirement that $\GG$ has two distinct elements. It follows that $g^{-1}\RR = g\RR$ and therefore 
$g^2\RR = (g\RR)(g\RR)= (g\RR)(g^{-1}\RR)=(gg^{-1})\RR=\RR$. We have obtained the important result that $g^2$ is a real number while $g$ itself is not. We can then set $g^2 = \alpha$, where $\alpha$ is a real number, not necessarily positive. Generalized complex numbers are obtained by generating a set that is closed under addition and multiplication of elements within cosets. Closure under addition generates the set 
$$\{a + gb | a,b \in \RR\},$$ for which
addition and multiplication rules follow directly from the properties of $g$:
\begin{eqnarray}
(a_1 + g b_1)+ (a_2 + g b_2) &=& (a_1 + a_2) + g(b_1 + b_2), \\
(a_1 + g b_1) (a_2 + g b_2) & = &  (a_1  a_2 +\alpha  b_1 b_2) + g(a_1 b_2 + a_2 b_1) .
\label{eq:complex}
\end{eqnarray}
We recognize that for $\alpha = -1$, the set just generated is the field of complex numbers, $\mathbb{C}$.

As typically done, Eq.~\ref{eq:complex} can be written in matrix form, 
\begin{equation}
\left( 
a_1 a_2+\alpha b_1 b_2,   a_1 b_2 + a_2 b_1   
\right) 
=
(a_1, b_1)
\left( \begin{array}{cc}
a_2  &  b_2\\
\alpha b_2 &  a_2
\end{array} \right), 
\end{equation}
which gives the general matrix form for complex numbers, 
\begin{equation}
z = 
\left( \begin{array}{cc}
a  & b\\
\alpha b &  a
\end{array} \right) = a \mathbf{1} + b \mathbf{I},
\label{eq:matrix-2d}
\end{equation}
where $\mathbf{1}$ is the identity matrix and 
\begin{equation}
\mathbf{I} =  
\left( \begin{array}{cc}
0  & 1\\
\alpha  &  0
\end{array} \right).
\label{eq-gen-2d}
\end{equation}

Although $\alpha$ can assume any real value, it is generally sufficient to consider the particular values  $-1, 0$, and $1$. This is because we have $\alpha\RR = \RR$ for any real value of $\alpha$ rendering the overall (positive) scale irrelevant for most purposes. The choice $\alpha = -1$ gives the usual complex numbers, $\alpha = 1$ the split-complex (or hyperbolic numbers) and $\alpha = 0$ the dual complex numbers \cite{yaglom}.

As it is known, the number $z$ has an inverse if $\det(z)= a^2-\alpha b^2$ is non zero, in which case the inverse element is $(a - g b)/\det(z)$. Split complex numbers do not form a field since multiplicative inverses do not exist for $a = \pm b$. Similarly, dual complex numbers do not have inverses for $a = 0$ (pure imaginary numbers). Also note that the transpose of the matrix in Eq.~\ref{eq:matrix-2d} is an equivalent representation of a complex number. In general, the matrix representations that we obtain below are determined up to a transposition. 

\section{Case $n = 3$} 

In this case the coset group $\GG$ is a cyclic group of order 3  with  $\GG = \{\RR, g\RR, g^{-1}\RR\}$.  The cosets $g\RR$ and $g^{-1}\RR$ are distinct, with $g^2\RR = g^{-1}\RR$, giving $g^{-2}\RR = g\RR$, 
and $g^3\RR = g^{-3}\RR = \RR$. By choosing elements $i \in g\RR$ and $j \in g^{-1}\RR$, we have the rules:
\begin{eqnarray}
ij & \in & \RR \nonumber \\
ji  & \in & \RR \nonumber \\
i^2  & \in & j\RR \nonumber \\
j^2 & \in & i\RR, 
\end{eqnarray}
from which it follows that $i^3 \in \RR$ and $j^3 \in \RR$. The above mean that there exist three real numbers $\alpha$, $\beta$, and $\gamma$ such that
\begin{eqnarray}
ij &= & ji = \alpha \nonumber \\
i^2  & = & \beta j \nonumber \\
j^2 & = & \gamma i.
\end{eqnarray}
Note that $i$ and $j$ commute because $g\RR$ and $g^{-1}\RR$ are inverse elements of each other in the coset group. (Commutativity can also be seen by setting $ij = \alpha$, $ji = \alpha'$ and by evaluating the product $iji$ which gives $\alpha i = i \alpha'$ and therefore $\alpha' = \alpha$.) 
Furthermore, parameters $\alpha$, $\beta$, and $\gamma$ are not independent. It can be shown that $\alpha = \beta \gamma$, for example by considering the product $i^2j = \beta j^2$ which gives $\alpha i = \beta \gamma i$. Table \ref{table-3d} lists possible assignments for parameters $\alpha$, $\beta$, and $\gamma$, with $\alpha = ij$ chosen to designate the signature of the number system.

\begin{table}[h]
\caption{Possible scaling assignments for 3D extensions of $\RR$. \newline}

\begin{tabular}{r|rr|r|r}
\hline \hline
$\alpha$ & $\beta$ & $\gamma$ & $i^3= \alpha \beta$ & $j^3 = \alpha \gamma$ \\
\hline
1 & 1 & 1 & 1 & 1 \\
 & $-1$ & $-1$ & $-1$ & $-1$ \\
\hline
0 & 1 & 0 & 0 & 0 \\
 & 0 & 0 & 0 & 0  \\
 & $-1$ & 0 & 0 & 0 \\
\hline
$-1$ & 1 & $-1$ & $-1$ & 1 \\
\hline \hline
\end{tabular}
\label{table-3d}
\end{table}

As shown in Table \ref{table-3d}, there are 6 distinct choices (of 0, -1, +1) for parameters $\alpha$, $\beta$ and $\gamma$. Some alternative assignments are equivalent and are not shown. For example, the assignment $\beta = 1$, $\gamma = -1$ is equivalent to $\beta = -1$, $\gamma = 1$ because $i$ and $j$ can be switched. 

Multiplying two generic 3D numbers of the form $a + bi + cj$, and setting the result in matrix form (as done above for complex numbers), gives the general matrix representation for 3D numbers as
\begin{equation}
z = 
\left( \begin{array}{ccc}
a  & b & c \\
\alpha c  & a & \beta b\\ 
\alpha b  & \gamma c & a 
\end{array} \right) 
= a \mathbf{1} + b \mathbf{I} + c \mathbf{J},
\label{eq:matrix-3d}
\end{equation}
where $\mathbf{1}$ is the identity matrix and
\begin{equation}
\mathbf{I} = 
\left( \begin{array}{ccc}
0  & 1 & 0 \\
0  & 0 & \beta \\ 
\alpha  & 0 & 0 
\end{array} \right), \, 
\mathbf{J} = 
\left( \begin{array}{ccc}
0  & 0 & 1 \\
\alpha  & 0 & 0 \\ 
0  & \gamma & 0 
\end{array} \right).
\label{eq-gen-3d}
\end{equation}

Since the general matrix form in Eq. ~\ref{eq:matrix-3d} is invariant under matrix multiplication, the form of the determinant, $det(z) = a^3 + \alpha \beta b^3 + \alpha \gamma c^3 - 3 \alpha ab c$ is invariant as well for all $\alpha = \beta \gamma$. 
The assignment $\alpha = \beta = \gamma = 1$ (a hyperbolic number system) is treated in detail in \cite{olariu-3d, olariu-book} where geometric interpretation as well as exponential and trigonometric forms are analyzed in addition to matrix representations.

\section{Case $n = 4$, cyclic}

The cyclic structure for $n = 4$ gives the following form for $\GG$, 
\begin{equation}
\GG = \{\RR, g\RR, g^{-1}\RR, g^2\RR = g^{-2}\RR\},
\end{equation}
in which $g^3\RR = g^{-1}\RR$, $g^{-3}\RR = g\RR$, $g^4\RR = g^{-4}\RR = \RR$.
With elements $i$, $j$, $k$ from $g\RR$, $g^{-1}\RR$, and $g^2\RR$, respectively, coset products within this group give the following multiplication rules:
\begin{eqnarray}
ij & = & ji = \alpha \nonumber \\
jk & = & \beta i \nonumber \\
ki & = & \gamma j \nonumber \\
i^2 & = & \delta k \nonumber \\
j^2 & = & \epsilon k \nonumber \\
k^2 & = & \varphi,
\end{eqnarray}
where Greek symbols indicate real coefficients. The six scaling parameters are not independent since it can be shown that 
\begin{eqnarray}
\varphi & = & \beta \gamma, \\
\beta \delta & = & \gamma \epsilon = \alpha.
\end{eqnarray}
This leaves four independent scales to be chosen as combinations of $-1$, $0$, and $1$. Table \ref{table-4dc} gives possible assignment values with $\alpha = ij$ chosen as the signature of the number system.

\newpage

\begin{table}[h]
\caption{Possible scaling assignments for 4D cyclic extensions of $\RR$. The last column indicates particular assignments considered in references \cite{olariu-4d, olariu-book}.\newline}
\begin{tabular}{r|rrrr|r|c}
\hline \hline
$\alpha$ & $\beta$ & $\gamma$ & $\delta$ & $\epsilon$ & $\varphi$ & \\ 
\hline
$1$ &  $1$ & $1$ & $1$ & $1$   & $1$ & Polar\\
   & $1$ & $-1$ & $1$ & $-1$ & $-1$ & \\
   & $-1$ & $1$ & $-1$ & $1$ & $-1$ & \\
   & $-1$ & $-1$ & $-1$ & $-1$ & $1$ & \\
\hline
0 & \multicolumn{4}{|c|}{$\beta \delta = 0$, $\gamma \epsilon = 0$} & $\beta \gamma$ & \\
\hline
$-1$ & $1$ & $1$ & $-1$ & $-1$ & $1$ & \\
   & $1$ & $-1$ & $-1$ & $1$ & $-1$ & \\
   & $-1$ & $1$ &  $1$ & $-1$ & $-1$ & Planar \\
   & $-1$ & $-1$ & $1$ & $1$ & $1$ & \\
\hline \hline
\end{tabular}
\label{table-4dc}
\end{table}
The general matrix form for cyclic 4D numbers is
\begin{equation}
z =
\left( \begin{array}{cccc}
a  & b & c & d\\
\alpha c  & a & \gamma d & \delta b\\ 
\alpha b  & \beta d & a& \epsilon c \\
\beta \gamma d & \beta c & \gamma b & a 
\end{array} \right) 
=  a \mathbf{1} + b \mathbf{I} + c \mathbf{J} + d \mathbf{K},
\label{eq:matrix-4dc}
\end{equation}
with $\mathbf{1}$ being the unit matrix and 
\begin{equation}
\mathbf{I} = 
\left( \begin{array}{cccc}
0  & 1 & 0 & 0\\
0  & 0 & 0 & \delta \\ 
\alpha & 0  & 0 & 0 \\
0 & 0 & \gamma & 0 
\end{array} \right),  
\mathbf{J} = 
\left( \begin{array}{cccc}
0  & 0 & 1 & 0\\
\alpha & 0  & 0 & 0 \\
0  & 0 & 0 & \epsilon \\ 
0 & \beta & 0 & 0 
\end{array} \right), 
\mathbf{K}  = 
\left( \begin{array}{cccc}
0  & 0 & 0 & 1\\
0 & 0 & \gamma & 0 \\ 
0 & \beta & 0 & 0  \\
\beta \gamma & 0  & 0 & 0
\end{array} \right).
\label{eq-gen-4dc}
\end{equation}

\section{Case $n = 4$, Klein group}

In this case, the coset group has the form $\mathbb{G} = \{\RR, g_1\RR=g_1^{-1}\RR, g_2\RR=g_2^{-1}\RR,  g_3\RR=g_3^{-1}\RR\}$, 
where $g_3 = g_1 g_2$. Choose elements $i$, $j$, $k$ belonging to $g_1\RR$, $g_2\RR$, and $g_3\RR$, respectively. With Greek symbols indicating scaling coefficients, coset products within this group give the following multiplication rules:
\begin{eqnarray}
i^2 & =& \alpha \nonumber \\
j^2 & = & \beta   \nonumber \\
k^2 & = & \gamma  \nonumber \\
jk & = & \alpha' i  \nonumber \\
ki & = & \beta' j  \nonumber \\
ij & = & \gamma' k  \nonumber \\
kj & = & \alpha'' i  \nonumber \\
ik & = & \beta'' j  \nonumber \\
ji & = & \gamma'' k.
\end{eqnarray}
The 9 scaling parameters are obviously not independent. Constraints are derived from products such as $ijk$, $(ij)(ji)$, and $(ij)(jk)$. We obtain:
\begin{eqnarray}
\alpha \beta & = & \gamma \gamma' \gamma'' \\
\beta \gamma & = & \alpha \alpha' \alpha'' \\
\gamma \alpha & = & \beta \beta' \beta'' \\
\alpha \alpha' &=& \beta \beta' = \gamma \gamma' = \alpha'' \beta'' \gamma'' \\
\alpha \alpha'' &=& \beta \beta'' = \gamma \gamma'' = \alpha' \beta' \gamma'.
\end{eqnarray}
An alternative way of expressing the constraints is
\begin{eqnarray}
\alpha'   \beta' & = & \alpha'' \beta'', \\
\alpha'  \beta'' & = & \alpha''  \beta' = \gamma, 
\end{eqnarray}
and all corresponding cyclic permutations. With notations $\rho = \alpha \beta \gamma$, $\rho' = \alpha' \beta' \gamma'$, and $\rho'' = \alpha'' \beta'' \gamma''$, it can be shown that $\rho \rho' \rho'' = \rho''^4 = \rho'^4 = \rho^2 \ge 0$. It means that the three products $\rho$, $\rho'$, and $\rho''$ are either all equal to zero or all different than zero. Note that $\rho'$ and $\rho''$ are products of basis elements: $\rho'' = ijk = jki = kij$ and $\rho' = ikj = kji = jik$. One convenient assignment procedure is to choose $\alpha$, $\beta$, $\gamma$, and $\rho''$. Then primed parameters are determined from $\alpha \alpha' = \beta \beta' = \gamma \gamma' = \rho''$ and double-primed parameters from $\alpha \alpha'' =  \beta \beta'' = \gamma \gamma'' = \alpha' \beta' \gamma' = \rho'$. In this parameterization, Hamilton's quaternions correspond to $\alpha = \beta = \gamma = \rho'' = -1$. This and other possible non-zero assignments are summarized in Table \ref{table-4dk}. The last column indicates particular assignments corresponding to common quaternionic number systems: quaternions (Q), split-quaternions (S), bicomplex numbers (B), and to 4D numbers analyzed in \cite{olariu-4d, olariu-book}, namely hyperbolic (H) and circular (C). The products $\rho$, $\rho'$, $\rho''$ can be considered as signatures for each assignment and can provide a quick consistency check. Also note that commutativity requires that corresponding primed and double-primed parameters are equal, which is possible if and only if $\rho > 0$. There are four commutative cases in Table \ref{table-4dk} including the bicomplex (B), hyperbolic (H), and circular (C) numbers, and four non-commutative cases including the quaternions (Q) and split-quaternions (S). 

\begin{table}[h]
\caption{Possible scaling assignments for 4D Klein extensions of $\RR$. Particular number systems are indicated in the last column:  hyperbolic (H), split-quaternions (S, called pseudoquaternions in Ref.~\cite{yaglom}), bicomplex (B), quaternions (Q), and circular (C).}

\begin{tabular}{rrr|r|rrr|r|rrr|r|l}
\hline \hline
$\alpha$ & $\beta$ & $\gamma$ & \mbox{ $\, \, \rho \, \, $ } &
$\alpha'$ & $\beta'$ & $\gamma'$ & \mbox{ $\, \, \rho' \, \, $ } &
$\alpha''$ & $\beta''$ & $\gamma''$ & \mbox{ $\, \, \rho'' \, \, $ } \\
\hline
1 & 1 & 1 & 1 & 1 & 1 & 1 & 1 & 1 & 1 & 1 & 1 &  H  \hspace{8pt}\\
$-1$ & 1 & 1 & $-1$ & $-1$ & 1 & 1 & $-1$ & 1 & $-1$ & $-1$ &  & S \\
$-1$ & $-1$ & 1 & 1 & $-1$ & $-1$ & 1 & 1 & $-1$ & $-1$ & 1 &   & B \\
$-1$ & $-1$ & $-1$ & $-1$ & $-1$ & $-1$ & $-1$ & $-1$ & 1 & 1 & 1 &  & \\
\hline
$-1$ & $-1$ & $-1$ & $-1$ & 1 & 1 & 1 & 1 & $-1$ & $-1$ & $-1$ & $-1$ & Q\\
1 & $-1$ & $-1$ & 1 & $-1$ & 1 & 1 & $-1$ & $-1$ & 1 & 1 &   & C \\
1 & 1 & $-1$ & $-1$ & $-1$ & $-1$ & 1 & 1 & 1 & 1 & $-1$ &  & \\
1 & 1 & 1 & 1 & $-1$ & $-1$ & $-1$ & $-1$ & $-1$ & $-1$ & $-1$ &  &  \\
\hline \hline
\end{tabular}
\label{table-4dk}
\end{table}
If we allow any of the parameters to be zero, then automatically we need to satisfy the constraint $\rho = \rho' = \rho'' =0$. These assignments give rise to degenerate quaternions (page 24 of \cite{yaglom}). Assignment options for zero-valued parameters are summarized in Table \ref{table-4d0} and particular assignments are indicated in the last column as follows: degenerate pseudoquaternions (DPQ), doubly degenerated quaternions (DDQ), and degenerate quaternions (DQ) \cite{yaglom}. 

\begin{table}[h]
\begin{threeparttable}[b]
\caption{Assignments for Klein 4D extensions of $\RR$ for $\rho = \rho' = \rho'' =0$. 
Particular number systems are indicated in the last column: degenerate pseudoquaternions (DPQ), doubly degenerated quaternions (DDQ), and degenerate quaternions (DQ).} 
\label{table-4d0}
\begin{tabular}{rrr|c|rrr|c|rrr|c|l}
\hline \hline
\mbox{ $\alpha$ }  & \mbox{ $\beta$ }  & \mbox{ $\gamma$ } 
& \mbox{ $\, \, \rho \, \, $ } &
\mbox{ $\alpha'$ } & \mbox{ $\beta'$ } & \mbox{ $\gamma'$ } 
& \mbox{ $\, \, \rho' \, \, $ } &
\mbox{ $\alpha''$ } & \mbox{ $\beta''$ }& \mbox{ $\gamma''$ }
& \mbox{ $\, \, \rho'' \, \, $ } & 
\\
\hline
1 & 0 & 0 & 0 & 0 & 1      & 1      & 0 & 0 & 1 & 1          & 0     & \\
  &    &    &    & 0 & $-1$ & 1      &    & 0 & 1 & $-1$      &  & DPQ \hspace{25pt} \\
  &    &    &    & 0 & $-1$ & $-1$ &    & 0 & $-1$ & $-1$ &  & \\
\hline
0 & 0 & 0 & 0 & 0 & 0 & $\pm 1$   & 0 & 0 & 0 & 0           & 0 & \\
   &    &    &   &  0  & 0 &  0             &   & 0 & 0 & $\pm 1$ &  & \\
   &    &    &   &  0  & 0 & $\pm 1$   &   & 0 & 0 & $\pm 1$  &  & DDQ\tnote{a}\\
   &    &    &   &  0  & 0 &   0           &    & 0 & 0 & 0            &  & \\
\hline
$-1$ & 0 & 0 & 0 & 0  & 1      & 1      & 0 & 0 & $-1$ & $-1$ & 0 & DQ\\
       &    &    &    & 0  & $-1$ & 1      &    & 0 & $-1$ & 1      &  & \\
       &    &    &    & 0  & $-1$ & $-1$ &    & 0 & 1 & 1           &  &  \\
\hline \hline
\end{tabular}
\begin{tablenotes}
\item[*] The $\pm$ signs do not need to correspond; all four possibilities are allowed. The DDQ assignment in \cite{yaglom} corresponds to the non-commutative case $\gamma' = -\gamma'' = 1$.
\end{tablenotes}
\end{threeparttable}
\end{table}

The general matrix form for Klein 4D numbers is
\begin{equation}
z  = 
\left( \begin{array}{cccc}
a  & b & c & d\\
\alpha b  & a & \beta'' d & \gamma' c\\ 
\beta c  & \alpha' d & a& \gamma'' b \\
\gamma d & \alpha'' c & \beta' b & a 
\end{array} \right) =  a \mathbf{1} + b \mathbf{I} + c \mathbf{J} + d \mathbf{K},
\label{eq:matrix-4dk}
\end{equation}
with
\begin{equation}
\mathbf{I} = 
\left( \begin{array}{cccc}
0  & 1 & 0 & 0\\
\alpha  & 0 & 0 & 0 \\ 
0 & 0  & 0 & \gamma'' \\
0 & 0 & \beta' & 0 
\end{array} \right), \, 
\mathbf{J} = 
\left( \begin{array}{cccc}
0  & 0 & 1 & 0\\
0  & 0 & 0 & \gamma' \\ 
\beta & 0  & 0 & 0 \\
0 & \alpha'' & 0 & 0 
\end{array} \right), \, 
\mathbf{K}
\left( \begin{array}{cccc}
0  & 0 & 0 & 1\\
0 & 0 & \beta'' & 0 \\ 
0 & \alpha' & 0 & 0  \\
\gamma & 0  & 0 & 0 
\end{array} \right).  
\label{eq-gen-4dk}
\end{equation}

Note that matrices corresponding to cyclic and Klein coset groups in Eq.~\ref{eq:matrix-4dc} and Eq.~\ref{eq:matrix-4dk} have different forms since the two groups have different structures. Therefore number systems from the two categories will have significantly different properties. For example, references \cite{olariu-4d, olariu-book} compare and give a detailed analysis of four different kinds of commutative 4D numbers called polar, planar, hyperbolic, and circular. The first two belong to cyclic cosets and the last two correspond to Klein cosets, as indicated in Tables \ref{table-4dc} and \ref{table-4dk}. 

\section{2D Extensions over the field of complex numbers}

The multiplication rules obtained above apply to extensions of complex numbers as well. In this case, coefficients and multiplication parameters can assume complex values. As an example, the coset group $\mathbb{G} = \{\mathbb{C}, g\mathbb{C}\}$ can be used to generate bicomplex numbers. (This is a particular case of Cayley-Dickson construction of algebras, also called Dickson doubling, e.g. in \cite{conway-smith}.) The simplest construction using Eq.~\ref{eq:matrix-2d} is the block matrix,
\begin{equation}
z = 
\left( \begin{array}{cc}
A  & B \\
\alpha B &  A
\end{array} \right), 
\end{equation}
where $A$, and $B$, and $\alpha$ are now complex numbers. Representing $A$ and $B$ as 2x2 real matrices,  we obtain
\begin{equation}
z = 
\left( \begin{array}{cc|cc}
a & b  & c & d \\
\alpha b & a & \alpha d & c \\
\hline
\alpha c  & \alpha d & a & b \\
\alpha^2 d & \alpha c & \alpha b & a
\end{array} \right) 
=  a \mathbf{1} + b \mathbf{I} + c \mathbf{J} + d \mathbf{K},
\end{equation}
with
\begin{equation}
\mathbf{I} = 
\left( \begin{array}{cccc}
0  & 1 & 0 & 0\\
\alpha  & 0 & 0 & 0 \\ 
0 & 0  & 0 & 1 \\
0 & 0 & \alpha & 0 
\end{array} \right), \, 
\mathbf{J} = 
\left( \begin{array}{cccc}
0  & 0 & 1 & 0\\
0  & 0 & 0 & 1 \\ 
\alpha & 0  & 0 & 0 \\
0 & \alpha & 0 & 0 
\end{array} \right), \, 
\mathbf{K}
\left( \begin{array}{cccc}
0  & 0 & 0 & 1\\
0 & 0 & \alpha & 0 \\ 
0 & \alpha & 0 & 0  \\
\alpha^2 & 0  & 0 & 0 
\end{array} \right).  
\end{equation}

Comparison with the general matrix form of Klein 4D numbers in Eq. ~\ref{eq-gen-4dk} gives the assignments shown in Table \ref{table-corresp}.
\newpage

\begin{table}[h]
\caption{Correspondence between 4D Klein extension over reals and 2D extension over complex. The last three rows show particular real assignments corresponding to $\alpha = -1$ (bicomplex, B), $\alpha = 0$ (dual, D), and $\alpha = 1$ (hyperbolic, H).\newline}

\begin{tabular}{rrr|c|rrr|c|rrr|c|l}
\hline \hline
\mbox{ $\alpha$ }  & \mbox{ $\beta$ }  & \mbox{ $\gamma$ } 
& \mbox{ $\, \, \rho \, \, $ } &
\mbox{ $\alpha'$ } & \mbox{ $\beta'$ } & \mbox{ $\gamma'$ } 
& \mbox{ $\, \, \rho' \, \, $ } &
\mbox{ $\alpha''$ } & \mbox{ $\beta''$ }& \mbox{ $\gamma''$ }
& \mbox{ $\, \, \rho'' \, \, $ } & \\
\hline
$\alpha$ & $\alpha$ & $\alpha^2$ & $\alpha^4$ & 
$\alpha$ & $\alpha$ & $1$ & $\alpha^2$  & 
$\alpha$ & $\alpha$ & $1$ & $\alpha^2$ & \\
\hline
-1 & -1 & 1 & 1 & -1 & -1 & 1 & 1 & -1 & -1 & 1 & 1 & B \hspace{8pt} \\  
0  & 0 & 0 & 0 & 0 & 0 & 1 & 0 & 0 & 0 & 1 & 0 & D    \\
1 & 1 & 1 & 1 & 1 & 1 & 1 & 1 & 1 & 1 & 1 & 1 & H    \\
\hline \hline
\end{tabular}
\label{table-corresp}
\end{table}

\section{Conclusions}

The coset approach allows a systematic construction of multidimensional number systems for which the basis and multiplication rules follow from the structure of the coset group.  Constraints on multiplication parameters are obtained from the associativity property of coset product. Although assignments to $-1$, $0$, and $1$ appear most natural, there is no restriction to these three values as long as multiplication constraints are satisfied. For example, the parameter $\alpha$ for 2D numbers can take any real value except that for $\alpha \ge 0$, a subset of numbers exist with no multiplicative inverses. Higher dimensional systems can be generated in a similar way by considering the appropriate coset group. For $n \ge 5$, the number of possible group structures can be large, except if $n$ is a prime number in which case  the coset group is necessarily cyclic. 

\section*{Acknowledgments}

The author thanks Drs. Kashyap Vasavada and Prashant Srinivasan for discussions and reading of the manuscript.

\section*{APPENDIX} 

\subsection*{1. Introduction to groups and subgroups}
This brief elementary introduction is based on standard textbook abstract algebra (\textit{e. g.} \cite{hungerford, deskins-algebra, ash-algebra, pinter-algebra}).
 
A \textbf{group} is a set of elements $G= \{a, b, c, ...\}$ which is closed under an operation between elements. Generically, this group operation is called product, and closure means that the product of any two elements in the group is itself an element of the group. In addition, the set must include an identity element and each element in the set must have an inverse. Formally, a group $G$ has the following four properties:

\textbf{1. Closure}: $ab \in G$, for all $a,b \in G$.

\textbf{2. Associativity}: $(ab)c = a(bc)$ for all $a,b,c \in G$.

\textbf{3. Identity element}: there exists $e \in G$ such that $ae = ea = a$ for all $a \in G$.

\textbf{4. Inverse elements}: for each $a \in G$, there exists and element $a' \in G$ such that
$aa' = a'a = e$. The element $a'$ is the inverse element of $a$ and it is usually written as $a^{-1}$.

A \textbf{subgroup} is a subset $N$ of a group such that $N$ is closed under the group operation. Briefly stated, a subgroup is a closed subset of a group. 

\subsection*{2. Introduction to cosets}
Given a group $G$, a subgroup $N$ of $G$  is called normal if it is invariant under conjugation, defined as
\begin{equation}
gNg^{-1} = N,\forall g \in G.
\end{equation}
Note that this kind of transformation does not require that each element in $N$ is preserved but only that the set $N$ as a whole is preserved.  However, if $g$ commutes with any element $n \in N$, then $n$ is invariant under conjugation since
$gng^{-1} = gg^{-1}n = n$.

A normal subgroup $N$ can be used to generate \textbf{cosets} through the action of outside elements $g \notin N$. Note that the new set is distinct from $N$ if and only if $g \notin N$, otherwise if $g \in N$ then $gN = N$ due to closure. The important property of cosets is that there is no partial overlap between them. An operation called coset product can be defined on the set of cosets namely
\begin{equation}
(g_1N)(g_2N) = (g_1 g_2)N,
\end{equation}
where $g_1$ and $g_2$ are elements in $G$. The set of all cosets generated by a normal subgroup $N$ forms a group under the coset product. In this manuscript we look for finite coset groups that can be generated by the set of real numbers. 


%

\end{document}